\documentclass[12pt]{article}
\usepackage{amssymb}

\title{One-parameter families of functions in the Pick class}
\author{Waldemar Pa\l uba\thanks{Partially
                         supported by a KBN grant no.
                         2 PO3A 041 15. }\\
                         \small Institute of Mathematics\\
                         \small Warsaw University\\
                         \small Banacha 2\\
                         \small 02-097 Warsaw, Poland\\
                         \small e-mail: paluba@mimuw.edu.pl
                         }
\date{}

\begin{document}

\newtheorem{df}{Definition}
\newtheorem{pr}{Proposition}
\newtheorem{conj}{Conjecture}
\newtheorem{theo}{Theorem}
\newtheorem{rem}{Remark}
\newtheorem{lem}{Lemma}
\newcommand{\vp}{\varphi}

\maketitle

\begin{abstract} In the one-parameter family of power-law maps of the form
$f_a(x)=-|x|^{\alpha}+a,$ $\alpha
>1,$ we give examples of mutually related dynamically determined
quantities, depending on the parameter $a$, such that one is a
Pick function of the following one. These Pick functions are
extendable by reflection through the $(1,+\infty)$ half-axis and
have completely monotone derivatives there.
\end{abstract}


\newpage

\section{Introduction}

The study of the monotonicity of the dynamics in a one parameter family
of power-law maps $f_a = -|x|^{\alpha} + a \ $ with generally non-integer
$\alpha > 1 \ $ and a real parameter $a \,$ has a long history.
This was for that family where the first results pertaining to the
uniqueness of some sort of dynamical behavior (like the Feigenbaum dynamics,
for instance) were established, though only for $\alpha \, $'s close enough
to $1 \, $, cf. \cite{CE}.
The case of special interest $\ \alpha = 2 \ $ was later successfully
solved with use of quasi-conformal technics.

Roughly speaking, the results for $\alpha$ close to 1 were
obtained by perturbations of the well understood linear toy-model.
The perturbation was getting out of control with increasing
$\alpha \,$ and the idea that the behavior for larger $\alpha
\,$'s has its roots in the case where there is nearly no
nonlinearity present did not enjoy wide recognition since then. We
hope this research notes can help revive an interest in that
approach to the problem, although the examples described here,
that we were able to solve completely, are dynamically very simple
ones.

Anyway, as far as we know, what we give here is the first example
of carrying the information obtained for $\alpha$ close to  $1$ up
through all greater than $1$ exponents, without loss of control of
the influence of the perturbation made.

Let us briefly describe the strategy employed. We use arguments
from complex variable as well as
 geometric properties of the
root mapping $z \mapsto z^{\frac{1}{\alpha}} \, $. The former
involves properties of the so called Pick functions i.e. the
mappings that take the upper half-plane into itself. Whenever such
a mapping is extendable by reflection throughout a half-axis and
its derivative vanishes at $+\infty \, $ than this derivative is
in the class of completely monotone functions. With some abuse of
the language we may say that in what follows we study dynamically
determined functions in the parameter space and their complex
extensions. The actual parameter is replaced with some related
quantity in $[1 \, , +\infty) \, $. To prove the monotonicity of
the derivative of such a function we show the function itself is
in the appropriate Pick class. To see the latter we examine
geometric properties of the root maps, making use of the
contraction of the Poincar\'e  neighborhoods about the rays (cf.
the idea stated in \cite{Sul} in the proof of the sector theorem),
as well as of the property of taking rays onto rays. In the course
of the proof we show that some regions are transformed in a
desired way, eventually leading to the Pick property of our maps
in the extended parameter space. Staying in the completely
monotone derivative class while perturbing $\alpha \, $ provides
for the control of the derivative and implies the existence of a
strictly monotone in parameter, dynamically determined function
for the perturbed $\alpha \,$ and eventually for all large
exponents, because after the perturbation we keep the derivative
of the inverse function uniformly bounded away from $0$. The
prospect of using this kind of technic  for more intricate
situations depend upon more complete understanding of the geometry
of the root maps at the infinitesimal level, for small
wedge-shaped neighborhoods. This is a work in progress.

Finally, let us relate the results of this paper to what we
previously knew about monotone behavior of dynamically determined
quantities in the power-law family. In our work \cite{P} we have
found an infinite sequence of monotone functions of the parameter.
That sequence had been determined by the whole post-critical orbit
rather than by just first few steps as it is the case in Theorems
1 and 2 below. The results of \cite{P} were achieved by purely
real variable means and functions in question differed from those
investigated in here. Instead, in the current work the information
derived from the complex plane extensions is much stronger.

\section{Complex variable prerequisites and geometry of the root map.}

Consider a complex variable function $\varphi$, analytic
in the upper half-plane $\bf H $. We denote by $P_{(1,+\infty)}$
the class of those Pick functions (ie. with  positive imaginary
part, cf. \cite{Don}) which can be continued by reflection
throughout the half-axis $(1, +\infty)$.

Assume that $\varphi\in P_{(1,+\infty)}$ and that $\varphi$ fixes $1$,
ie. it can be continuously defined along the reals
for $z=1$ and $\varphi(1)=1$.

We shall denote by $P$ the class of functions satisfying the above
and call it the class of Pick Argument Lessening functions, to
distinguish them from better recognized argument decreasing maps.
The reason for this name explains in Proposition  \ref{prjeden}
below.

A function in this class will simply be refered to as $a \  pal$.

\begin{rem}\label{remjeden}
Every function in the Pick class $P_{(a,+\infty)}$, i.e.
continuable throughout a half-axis
$(a,+\infty)$, that has a finite limit of the derivative at infinity
 has a completely monotone derivative on that half-axis
{\rm (cf. \cite{Don})}.
\end{rem}

The above property shall be the key to applying the class of $pals$
to the study of the problems of monotonicity of dynamical
behavior in a real valued parameter, like  in the
family $-|x|^{\alpha}+a,\ \alpha>1$.

Our main point throughout this work will be to show that inverse $pal$
functions show up naturally in the one-parameter power-law family
as some dynamically defined quantities with the parameter as the variable.

To that goal we will need an appropriate choice of variables, as
well as some facts concerning real variable properties of the
power-law mapping $x\mapsto |x|^{\alpha}$. This we postpone for
the moment, and instead  we focus here on preparatory complex
variable arguments.

The following proposition points out to
a relevant property of the maps in our class.

\begin{pr}\label{prjeden}
If $\varphi\in P$ then for every $z\in \bf H$ $\mbox{Arg} (\varphi
(z) - 1) \leq \mbox{Arg} (z - 1)$.
\end{pr}

\subparagraph*{Proof.}

For a number $z$ in the upper half-plane $Arg (z-1) = \alpha$
means that $z$ lies on the line crossing reals at $z=1$ at the
angle $\alpha$. By Schwarz-Pick  lemma a region of the upper
half-plane cut out by a circle tangent on the right hand side to
that line at the point $z=1$ is mapped by $\varphi$ into another
such a region. Letting the radius of the cutting circle tend to
infinity we get the statement of the proposition (cf. \cite{Sul},
proof of the sector theorem).

$\Box$

We now proceed to the study of the simpliest dynamical setting
where there appear  the Pick argument lessening functions.

Let $\alpha >1$ be a fixed real number. Consider the family of power-law
mappings $f_a = -|x|^{\alpha} + a \, , \, 0\leq a <1$.
The quantity $(1 - a^{\alpha -1})^{-1}$ is monotone in $a$ and so
this quantity itself may well serve as a parameter running the
$[1,+\infty)$ half-axis, rather than the interval $[0,1)$.

Notice that this number, which we shall denote by $p_2$,
is precisely the exponent of the length of the interval
$\left({f_a}^2(0) \, , \,  f_a(0)\right)$ measured in the Poincar\'e
metric on the positive half-axis.

Before moving on to the geometric part of this section we introduce some
more notation.

The interval $\left({f_a}^{n-1}(0), {f_a}^n(0)\right)$,
$n=3,4,\ldots$ is contained within the interval
$({f_a}^{n-2}(0) , {f_a}^{n-3}(0))$.

Denote by $p_n$ the exponent of the Poincar\'e length of the former
interval within the latter, i.e. the value of the cross-ratio

$$
\frac{{f_a}^{n-3}(0) - {f_a}^n(0)}{{f_a}^{n-3}(0) - {f_a}^{n-1}(0)}
\, \cdot \,
\frac{{f_a}^{n-1}(0) - {f_a}^{n-2}(0)}{{f_a}^n(0) - {f_a}^{n-2}(0)}
$$

For a point $\varphi\in{\bf H}$ we denote by $D_{\vp}$
the closed disc bounded by the circle $C_{\vp}$
passing through the points $0,1,\vp.$
The disc $D_{\vp}$ is the image of the upper half-plane
under the linear fractional transformation
$\ell_{\vp}:$
\[
\ell_{\vp}\!:t\mapsto\frac{t}{1+\frac{t-1}{\vp}}.
\]
Points $t$ in the upper half-plane satisfying the condtition
\[\mbox{Arg}(t-1)\geq \mbox{Arg}(\vp -1)\]
are mapped into the region between the circle $C_{\vp}$
and the circle passing through the points 1, $\vp,$
tangent to the ray connecting 0 and $\vp.$
We denote this region by $D_{\vp}^-.$
Here  we give a geometric lemma.

\begin{lem}\label{lematdwa}
For every $\alpha \geq 1$ the root mapping $r(z)= z^{\frac{1}{\alpha}}$
maps the region $D_{\vp}^-$ into the disc $D_{r(\vp)}.$
\end{lem}

\subparagraph*{Proof:}
We will jointly use Proposition 1 and the fact that the root
map takes rays originating at 0 onto rays.

The lower part of $D_{\vp}$ (i.e. the part lying
below the real axis) is mapped into itself, because the root
map is in the $P_{(0,+\infty)}$ class and fixes 0 and 1.
For the same reason the point $\vp$ is mapped onto a point
$r(\vp)$ in the upper part of $D_{\vp}.$

Thus the lower part of $D_{r(\vp)}$ contains the lower part
of $D_{\vp}$ and so it contains the whole image of the lower
part of $D_{\vp}.$
To do with the upper part we have to use the fact that the
root map is in $P.$
By Proposition 1 the point $r(\vp)$ lies to the right-hand
side of the ray originating at 1 and passing through $\vp.$
For any point $t\in {\bf H}$ the circle $C_t$ is intersected by the
ray connecting 0 and $t$ at the angle equal to $\mbox{Arg}(t-1)$,
and so this angle gets decreased when we replace $t$ by
$r(t).$

Therefore the part of the disc $D_{r(\vp)}$ above the ray
connecting 0 and $r(\vp)$ contains within itself the arc of
the circle that joints 0 and $r(\vp)$ and intersects
this ray at the angle $\mbox{Arg}(\vp -1).$  But the image
of the part of the disc $D_{\vp}$ lying above the ray
passing through 0 and  $\vp$  is totally contained
between the previously described arc and the ray connecting
0 and $r(\vp).$
Here we used the fact that the root map takes rays originating at 0 onto
rays and Proposition 1 together.

The image of the interval $(0,\vp)$ is a concave curve.
This completes the proof of the lemma also for the upper part
of the disc.
$\Box$

\begin{rem}
The statements of both Proposition 1 and Lemma 1 can be obtained
by direct computation in polar coordinates, employing
differentiating in $\alpha.$ and Jensen inequality.
\end{rem}
However, those computations  are rather tedious, whereas our
geometrical argument is much simpler.

\section{The Pick functions in the parameter space}

For a linear moder of our dynamical system we immediately
find a direct formula relating $p_2$ to $p_3.$

 \[
 p_2=\frac{1}{2}+\left(p_3-\frac{3}{4}\right)^{\frac{1}{2}}.\]

This function has a holomorphic extension
to the upper half-plane and this extension
is in the Pick class $P.$

If the exponent of the power law is greater than 1 we see that
the relation between  $p_2$ and $p_3$ is subject to the
following functional equation:
\begin{equation}\label{p2}
p_2=\left( 1+\frac{p_3-1}{p_2}\right)^{\frac{1}{\alpha}}.
\end{equation}
It is clear that for $\alpha$ close to 1
this equation has a solution which is a real analytic function
on $(1,+\infty),$ and that $p_2(1)=1.$ We claim this solution is
in $P.$

To see this we will consider a map $\vp(z)$ in $P$ and an operator
\[
R(\vp)(z)=\left( 1+\frac{z-1}{\vp(z)}\right)^{\frac{1}{\alpha}}.\]
Notice that this operator takes the class $P$ into itself.
Actually, if $\vp\in P$ and $z$ is a real greater that 1, so is
$R(\vp)(z).$

For $z\in{\bf H}$
\[
\mbox{Arg}(z-1) \geq \mbox{Arg}(\vp-1) > \mbox{Arg}(\vp) >0,\]
thus $\frac{z-1}{\vp}\in{\bf H},$
so $(R(\vp)(z))^{\alpha}\in {\bf H}$ and so does $R(\vp)(z).$
It follows that $R(\vp)$ is well defined for all $z\in {\bf H},$
maps the upper half-plane into itself, is extendable
by reflection throughout $(1,+\infty)$ and of course $R(\vp)(1)=1.$
Therefore $R(\vp)\in P.$

We yet have to investigate the convergence of the series of averages
\begin{equation}\label{star}
\frac{1}{n}\sum_{i=0}^{n-1} R^i(\vp)(z).
\end{equation}

If we begin with $\vp= Id,$ then all summands in (\ref{star})
will be uniformly bounded on compact subintervals of $(1,+\infty).$
Thus, due to a theorem about convergence of Pick functions
(cf. Chapter 2.4 in \cite{Don}), we can pick up a sequence of
(\ref{star}), convergent in the Pick class. The limit function of
this subsequence solves the functional equation
(\ref{p2}), but the solution $p_2(p_3)$ of
(\ref{p2}) is a uniquely defined analytic function on
$(1,+\infty)$, so it must belong to the Pick class itself.

At infinity, the derivative of $p_2$ tends to 0, so the function
$p_2(p_3)$ has a completely monotone derivative, according
to Remark \ref{remjeden}. Since the limit value of this derivative
at  the point 1 is smaller than 1, the inverse function $p_3(p_2)$
has the derivative greater than 1 all the time. Small increase
of the exponent $\alpha$ keeps the derivative of $p_3(p_2)$
positive all the time, so the inverse is, by the implicit function
theorem, a well defined analytic solution of the equation
(\ref{p2}) on $(1,+\infty).$
As we have seen, it again  has to be in the Pick class, with its derivative
completely monotone. So again its inverse function, that {\em a priori}
 merely had a
derivative positive, must actually have a derivative greater than 1 for
all its arguments. This shows that the property of bounded away from 1
derivative is kept while increasing $\alpha,$  so $\alpha$ can be
increased indefinitely and we see that for all $\alpha'$s larger than 1
the solution of the functional equation (\ref{p2}) is a well defined
analytic function in $P.$ This accounts for the proof of the
following:

\begin{theo}\label{thjeden}
For all exponents $\alpha >1$ the mapping $p_2(p_3)$ is in the
Pick Argument  Lessening class. \hfill$\Box$
\end{theo}

With some change in the proof, we could have delt with all
$\alpha'$s greater than $1$ simultaneously. We choose the above
argument to make the proof of Theorem 1 a preparatory step for
Theorem 2 below.

Now we proceed to the proof of the harder part.

\begin{theo}\label{thdwa}
For all exponents $\alpha >1$ the mapping $p_3(p_4)$
is in the Pick Argument Lessening class.
\end{theo}

\subparagraph*{Proof:}
The outline of the proof is much like in the previous theorem.
However the functional equation we deal with now is somewhat
different, and to see that the implied operator on Pick
functions has the range in the Pick class also, we shall need the
geometric properties of the root map stated in Section 2.

For the linear model of the dynamics, the function $p_3(p_4)$
is simply the identity. For a fixed $\alpha$ close to 1 the map in question
is a small perturbation of identity. Due to the dynamics relating
$p_4$ to $p_3$ it is elementary to see that the following
functional equation must be satisfied:

\begin{equation}\label{star2}
p_3(p_4)=\ell^{-1}_{p_2(p_3)}\left(\left(\frac{p_4}{1+\frac{p_4-1}
{p_2(p_3))^{\alpha}}}\right)^{\frac{1}{\alpha}}\right).
\end{equation}

From the proof of Theorem 1 we know not only that $p_2(p_3)\in{\bf H}$
for $p_3\in{\bf H},$ but also $(p_2(p_3))^{\alpha}\in{\bf H},$ and
$\mbox{Arg}((p_2(p_3))^{\alpha}-1)\leq \mbox{Arg}(p_3-1).$

Consider a number $w\in{\bf H}$ and suppose that $z=z(w)$
is a function in $P.$
Let $\psi(z)=p_2(z)^{\alpha}$ be the Pick function
defined for that fixed $\alpha$ in the proof of Theorem 1.

We define an operator on Pick functions:
\[
{\cal R}(z)(w)=\ell^{-1}_{\psi(z)^{\frac{1}{\alpha}}}
\left(
\left(
\frac{w}{1+\frac{w-1}{\psi(z)}}\right)^{\frac{1}{\alpha}}
\right).
\]
As before, we will see that its range is in $P.$ For a function
$z(w)\in P$ we have
\[
\mbox{Arg}(\psi(z)-1)\leq \mbox{Arg}(z-1)\leq \mbox{Arg}(w-1).
\]
The linear fractional map
$w\mapsto \frac{w\psi(z)}{\psi(z)+w-1}$ takes the number $w$ into
the region $D_{\psi(z)}^-,$ which, by Lemma 1, is mapped into the
disc $D_{p_2(z)},$ and after applying $\ell_{p_2(z)}^{-1}$
back into the upper half-plane.

So the function ${\cal R}(z)$ is in the Pick Argument Lessening class and
consequently $\mbox{Arg}(\psi({\cal R}(z))-1)\leq \mbox{Arg}(w-1),$
which makes the argument iterative.

Having dealt with the problem of remaining in the Pick class,
we are now in a position to mimic the further steps of the proof
of Theorem 1. We consider the sequence
\[
\frac{1}{n}\sum_{i=0}^{n-1}{\cal R}^i(z)(w),\]
starting with $z(w)=w.$  Its consecutive summands
 stay uniformly
bounded  on compact subsets of
$(1,+\infty).$ This again allows for eliciting
a convergent subsequence, which gives that a solution to (\ref{star2})
is in the Pick class. Again, looking at the inverse function, $p_4(p_3)$
on $(1,+\infty),$ we argue that its derivative is bounded away from 1,
so we can perturb it by increasing the exponent $\alpha,$ still
keeping the derivative bounded away from 0.
Again, looking at the inverse of a solution to (\ref{star2}) (with that
new $\alpha)$  we see that it has a non-vanishing derivative.
This solution has a completely monotone derivative, starting with
a value smaller that 1 and vanishing at $+\infty$, so its inverse actually
had a derivative bounded away from 1 rather than merely from 0.
This way, lifting up $\alpha $ indefinitely, we get a solution in the Pick class
for an arbitrary $\alpha >1.$
$\Box$

We clearly see that this kind of argument, that starts with
understandable behavior for $\alpha$ sufficiently close to 1,
and increasing the exponent while keeping the derivative bounded
away from 1 can be applied whenever we know that the derivative remains a
monotone function. In Theorems 1 and 2 this was provided for
by proving much stronger a statement, namely that it was a
completely monotone function.

The main obstacle for iterative use of the ideas from the proofs
of Theorems 1 and 2 is that the functional equation we get for
further steps requires understanding of the geometry of the root
maps acting on circles that pass through 1, but no longer pass
through 0. In this case we cannot use the argument of Lemma 1,
because the rays originating at 1 no longer are mapped onto rays.
So far we have not found a good remedy for this difficulty.

Actually, one might expect that for an arbitrary dynamics in the
power-law family the change of behavior while varying the
parameter is linked to some underlying completely monotone
function.

\end{document}